# A Less Noise-Sensitive SDP Relaxation in Wireless Sensor Network Localization

Pouya Mollaebrahim Ghari, Reza Shahbazian, Seyed Ali Ghorashi

***Abstract*—There are variety of methods to solve the localization problem and among them semi-definite programming based methods have shown great performance in both complexity and accuracy aspects. In this paper, we introduce a class of less noise-sensitive relaxation to reduce the complexity of the SDP-based methods. We apply our relaxation to Edge-based Semi-definite Programming method (ESDP) and the resulted model called PESDP. Simulation results confirm that our proposed PESDP method is less noise-sensitive and faster compared to the original ESDP.**

## I. Introduction

Nowadays, wireless sensor networks are considered to provide reliable solutions to a wide variety of applications including structural health monitoring, traffic control [1], industrial automation [2] and robotics [3]. We can obtain more purposeful data collected by a node only if we know its location. Therefore, the localization can be viewed as a necessity for wireless sensor networks.

The position of sensors can be determined by using a GPS system, but this could be expensive or an impossible solution in some cases [4]. However, the location of each node in a sensor network can be estimated based on the measurements of distances between neighboring nodes. In addition, there are a few nodes with known positions (called anchor) that can be used to solve the localization problem.

The localization problem can be described mathematically as follows. Here, we consider a two-dimensional (2D) localization problem whose extension to higher dimensions is straightforward. There are $n$ sensors with unknown locations and $m$ anchors whose locations are known as $a_1, \dots, a_m$. We define the Euclidean distance $d_{ij}$ for a pair of sensors $x_i$ and $x_j$, when the distance between them is less than the radio range. Similarly, for a sensor $x_j$ and anchor $a_k$, the Euclidean distance is denoted as $d_{kj}$. Therefore, we may write the localization problem as follows:

$$find \qquad X \in \mathcal{R}^{2\times n} \tag{1}$$

$$s.t. \, Y_{ii} - 2Y_{ij} + Y_{jj} = d_{ij}^2 \,, \forall (j,i) \in N_s$$

$$Y_{jj} - 2x_j^T a_k + \|a_k\|^2 = d_{jk}^2 \,, \forall (j,k) \in N_a$$

$$Y = X^T X \tag{1.d}$$

where $X = [x_1, \dots, x_n]$, $N_s = \{(j,i) \big| \|x_j - x_i\| < r\}$, $N_a = \{(j,k) \big| \|x_j - a_k\| < r\}$ and radio range is denoted by $r$. Using convex relaxation techniques is a very powerful approach to solve sensor network localization problems. Problem (1) is not a convex optimization problem. However, there are variety of relaxations that can transform it into a convex one. Semi-Definite Programming (SDP) relaxation which proposed in [5] is a powerful approach to solve the localization problem. Several methods have been proposed in order to enhance the accuracy of SDP [6]-[8]. Authors [9], [10] modified SDP in order to find a low rank solution, but these methods could not provide more accurate solution compared to SDP method introduced in [5]. Theoretical characteristics of SDP-based methods have been studied in [11], [12]. In such approaches, constraint (1.d) is relaxed to:

$$Y \succcurlyeq X^T X \quad \rightarrow \quad Z \succcurlyeq 0 \tag{2}$$

where $Z = \begin{pmatrix} I_2 & X^T \\ X & Y \end{pmatrix}$.

Edge-based Semi-Definite Programming (ESDP) relaxation is much faster than the original SDP and also comparable in terms of accuracy [13]. By applying ESDP relaxation to problem (1) we may write the localization problem as follows:

$$min \sum_{(j,i)\in N_s} (\alpha_{ij}^+ + \alpha_{ij}^-) + \sum_{(j,k)\in N_a} (\alpha_{jk}^+ + \alpha_{jk}^-) \tag{3}$$

$$s.t. \; diag(A_I^T Z A_I) = b_I$$

$$\begin{pmatrix} e_i - e_j \\ \mathbf{0} \end{pmatrix}^T Z \begin{pmatrix} e_i - e_j \\ \mathbf{0} \end{pmatrix} - \alpha_{ij}^+ + \alpha_{ij}^- = d_{ij}^2 \,, \forall (j,i) \in N_s \tag{3.c}$$

$$\begin{pmatrix} e_j \\ -a_k \end{pmatrix}^T Z \begin{pmatrix} e_j \\ -a_k \end{pmatrix} - \alpha_{jk}^+ + \alpha_{jk}^- = d_{jk}^2 \,, \forall (j,k) \in N_a \tag{3.d}$$

$$Z_{(1,2,i,j),(1,2,i,j)} \succcurlyeq 0\ , \forall (j,i) \in N_s$$

$$\alpha_{ij}^{+}, \alpha_{ij}^{-}, \alpha_{jk}^{+}, \alpha_{jk}^{-} \geq 0$$

$$i,j = 1,\dots,n\ , k = 1,\dots,m$$

Performance enhancement of ESDP method is more studied in [14]–[16].

In practical scenarios, measured distances are corrupted by noise and this can degrade the accuracy of the localization, especially when noise level is high. In this paper we modify ESDP relaxation in order to introduce a new less noise-sensitive relaxation and find a low rank solution. The rank minimization of matrix in (2) is usually done by means of objective function [9], [10]. In this paper we introduce a new method for rank minimization using the dual of the localization problem (3).

The remainder of the paper is organized as follows. Section II presents modified ESDP relaxation. In section III the numerical results are displayed and finally, section IV concludes the paper.

## II. Proposed convex relaxation

In this paper, we aim to modify ESDP relaxation in order to make it less noise-sensitive especially when noise level is high.

In the presence of noise, constraints (3.c) and (3.d) are rewritten as follows:

$$\begin{pmatrix} e_i - e_j \\ 0 \end{pmatrix}^T Z^{(n)} \begin{pmatrix} e_i - e_j \\ 0 \end{pmatrix} - \alpha_{ij}^{+} + \alpha_{ij}^{-} = d_{ij}^2 + u_{ij}\ , \forall (j,i) \in N_s \quad (4)$$

$$\begin{pmatrix} e_j \\ -a_k \end{pmatrix}^T Z^{(n)} \begin{pmatrix} e_j \\ -a_k \end{pmatrix} - \alpha_{jk}^{+} + \alpha_{jk}^{-} = d_{jk}^2 + v_{jk}\ , \forall (j,k) \in N_a$$

$$t_{ij} = 2n_{ij}d_{ij} + n_{ij}^2, \forall (j,i) \in N_s$$

$$v_{jk} = 2\delta_{jk}d_{jk} + \delta_{jk}^2, \forall (j,k) \in N_a$$

where additive noise associated with $d_{ij}$ and $d_{jk}$ are denoted by $n_{ij}$ and $\delta_{jk}$ , respectively and $Z^{(n)}$ denotes the noisy $Z$ matrix. When the measured distances are exact, we have:

$$\begin{pmatrix} e_i - e_j \\ 0 \end{pmatrix}^T Z^{(true)} \begin{pmatrix} e_i - e_j \\ 0 \end{pmatrix} = d_{ij}^2, \forall (j,i) \in N_s \quad (5)$$

$$\begin{pmatrix} e_j \\ -a_k \end{pmatrix}^T Z^{(true)} \begin{pmatrix} e_j \\ -a_k \end{pmatrix} = d_{jk}^2, \forall (j,k) \in N_a$$

From (4) and (5), we can conclude that the presence of noise can cause perturbation in $Z$ and consequently, the optimal value of (3) becomes larger and degrades the accuracy of the localization. Therefore, we can write:

$$Z^{(n)} = Z^{(true)} + \Delta \quad (6)$$

Now we consider dual problem of (3):

$$\max \sum_{(j,i)\in N_s} \omega_{ij}(d_{ij} + n_{ij})^2 + \sum_{(j,k)\in N_a} \omega_{jk}(d_{jk} + \delta_{jk})^2 + u_{11} + 2u_{12} + u_{22} \quad (7)$$

$$\text{s.t.} \sum_{(j,i)\in N_s} \omega_{ij}\left(0; e_i - e_j\right)^T \left(0; e_i - e_j\right) + \sum_{(j,k)\in N_a} \omega_{jk}(-a_k; e_j)^T(-a_k; e_j) + \begin{pmatrix} u_{11} + u_{12} & u_{12} & 0 \\ u_{12} & u_{22} + u_{12} & 0 \\ 0 & 0 & 0 \end{pmatrix} + \sum_{(j,i)\in N_s} S^{(i,j)} = 0$$

$$S_{(1,2,i,j),(1,2,i,j)}^{(i,j)} \succcurlyeq 0, \forall (j,i) \in N_s$$

$$S_{kl}^{(i,j)} = 0, \forall k \notin \{i,j\}\ or\ l \notin \{i,j\}$$

And we may write [17]:

$$\nabla_{\Delta_{(1,2,i,j),(1,2,i,j)}} p^*(0,0,0) = S_{(1,2,i,j),(1,2,i,j)}^{(i,j)^*} \quad (8)$$

The optimal value associated with (4) is denoted by $p^*(t, v, \Delta)$ and $S_{(1,2,i,j),(1,2,i,j)}^{(i,j)^*}$ is the optimal dual variable of (4). This means that if the absolute values of elements in $S_{(1,2,i,j),(1,2,i,j)}^{(i,j)^*}$ are decrease effectively, $p^*(t, v, \Delta)$ does not increase rapidly in the presence of noise. To do this, we modify ESDP method as follows:

$$min \sum_{(j,i)\in N_s} (\alpha_{ij}^{+} + \alpha_{ij}^{-}) + \sum_{(j,k)\in N_a} (\alpha_{jk}^{+} + \alpha_{jk}^{-}) \quad (9)$$

$$\text{s.t. } diag(A_I^T Z A_I) = b_I$$

$$\begin{pmatrix} e_i - e_j \\ 0 \end{pmatrix}^T Z \begin{pmatrix} e_i - e_j \\ 0 \end{pmatrix} - \alpha_{ij}^{+} + \alpha_{ij}^{-} = d_{ij}^2\ , \forall (j,i) \in N_s$$

$$\begin{pmatrix} e_j \\ -a_k \end{pmatrix}^T Z \begin{pmatrix} e_j \\ -a_k \end{pmatrix} - \alpha_{jk}^{+} + \alpha_{jk}^{-} = d_{jk}^2\ , \forall (j,k) \in N_a$$

$$Z_{(1,2,i,j),(1,2,i,j)} + P_{(1,2,i,j),(1,2,i,j)} \succcurlyeq 0\ , \forall (j,i) \in N_s \quad (9.e)$$

$$\alpha_{ij}^{+}, \alpha_{ij}^{-}, \alpha_{jk}^{+}, \alpha_{jk}^{-} \geq 0$$

$$i,j = 1,\dots,n\ , k = 1,\dots,m$$

By applying (9.e) to the optimization problem, all of the constraints in (7) remain unchanged. However, the objective function changes to the following:

$$\max \sum_{(j,i)\in N_s} \omega_{ij}(d_{ij} + n_{ij})^2 + \sum_{(j,k)\in N_a} \omega_{jk}(d_{jk} + \delta_{jk})^2 + u_{11} + 2u_{12} + u_{22} - \sum_{(j,i)\in N_s} tr(P_{(1,2,i,j),(1,2,i,j)} S_{(1,2,i,j),(1,2,i,j)}^{(i,j)}) \quad (10)$$

The optimal value of the problem is made robust to the perturbation in $Z^{(n)}$ by using relaxation in (9.e).

Now, we determine the perturbation matrix $P$ in order to find a low rank solution. Assume that $Z$ is a solution to (9) and $\{S^{(i,j)}\}$ is an optimal solution to the dual problem. Then, we have:

$$rank\left(Z_{(1,2,i,j),(1,2,i,j)}\right) \leq q \Leftrightarrow rank\left(S^{(i,j)}_{(1,2,i,j),(1,2,i,j)}\right) + 4 \leq 2q \quad (11)$$

The proof of (11) is given in Appendix A.

From (11) we can conclude that by minimizing the rank of $S^{(i,j)}_{(1,2,i,j),(1,2,i,j)}$, we can minimize the rank of $Z_{(1,2,i,j),(1,2,i,j)}$. It is known that if a matrix is symmetric and positive semi-definite we may minimize its trace as a convex approximation of its rank. Therefore, we aim to determine perturbation matrix $P_{(1,2,i,j),(1,2,i,j)}$ in order to minimize the rank of $S^{(i,j)}_{(1,2,i,j),(1,2,i,j)}$. Thus, the perturbation matrix is chosen as follows:

$$P_{(1,2,i,j),(1,2,i,j)} = p_{ij} I_4, \forall (j,i) \in N_s \quad (12)$$

Then we may rewrite (10) as follows:

$$\max \sum_{(j,i)\in N_s} \omega_{ij} (d_{ij} + n_{ij})^2 + \sum_{(j,k)\in N_a} \omega_{jk} (d_{jk} + \delta_{jk})^2 + u_{11} + 2u_{12} + u_{22} - \sum_{(j,i)\in N_s} p_{ij} tr(S^{(i,j)}_{(1,2,i,j),(1,2,i,j)}) \quad (13)$$

Therefore, by perturbation matrix in (12) we may minimize the rank of $S^{(i,j)}_{(1,2,i,j),(1,2,i,j)}$. This minimizes the rank of $Z_{(1,2,i,j),(1,2,i,j)}$ and consequently the computation complexity of the method reduces.

## III. SIMULATION RESULTS

In this section, several numerical comparisons for formulation (9) are reported. We evaluate the performance of PESDP in the presence of high level noise.

We consider two-dimensional (2D) localization problems and use benchmark test 10-500 which is available online at http://www.stanford.edu/~yyye/. In addition, we use MATLAB for simulations and SDPT3 solver in CVX software to perform our simulations [18]. We compute the position error for each network as follows:

$$\delta = \frac{1}{n}\sum_{i=1}^{n} \|\hat{x}_i - x_i\|$$

where the estimated location of the $i^{th}$ sensor is denoted by $\hat{x}_i$ and similarly, $x_i$ is a true position for this sensor. Therefore, we define the average position error as follows:

$$PE = \frac{1}{L}\sum_{l=1}^{L} \delta_l$$

where, L is the number of networks.

In figure 1, the effect of high noise level on the accuracy of ESDP, EML and our proposed method (PESDP) is studied. The distance measurements are corrupted by additive Gaussian noise. The maximum number of neighbors for each sensor is limited to 5. Two sensors are considered neighbors if corresponding distance exists. Networks consist of 300 sensors and 5 anchors and radio range is set to 0.2. 50 networks are simulated and perturbation matrix in (12) with $p_{ij} = 0.1$ is used for all relaxations. Figure 1 illustrates that our proposed PESDP obtains a better accuracy compared to EML and ESDP methods. In addition, as long as the standard deviation of noise increases, the difference between the accuracy of the proposed PESDP and other methods becomes larger. The perturbation matrix diminishes the effect of perturbation in constraints of the optimization model and as can be seen in figure 1, in the presence of high level noise, PESDP may obtain a better accuracy in comparison with ESDP.

In figure 2, we report the solution time of ESDP, EML and our proposed PESDP method by changing the number of sensors. The standard deviation of additive Gaussian noise is set to 0.1 and other properties are similar to prior simulation. As depicted in figure 2, the solution time of our proposed PESDP is less than the other methods and EML has a higher level of complexity compared with ESDP and our proposed PESDP. As can also be seen in figure 2, when the number of sensors increases, the difference between the solution time of PESDP and solution time of ESDP becomes larger. Therefore, simulation results confirm that the complexity of our proposed PESDP is less than the complexity of ESDP.

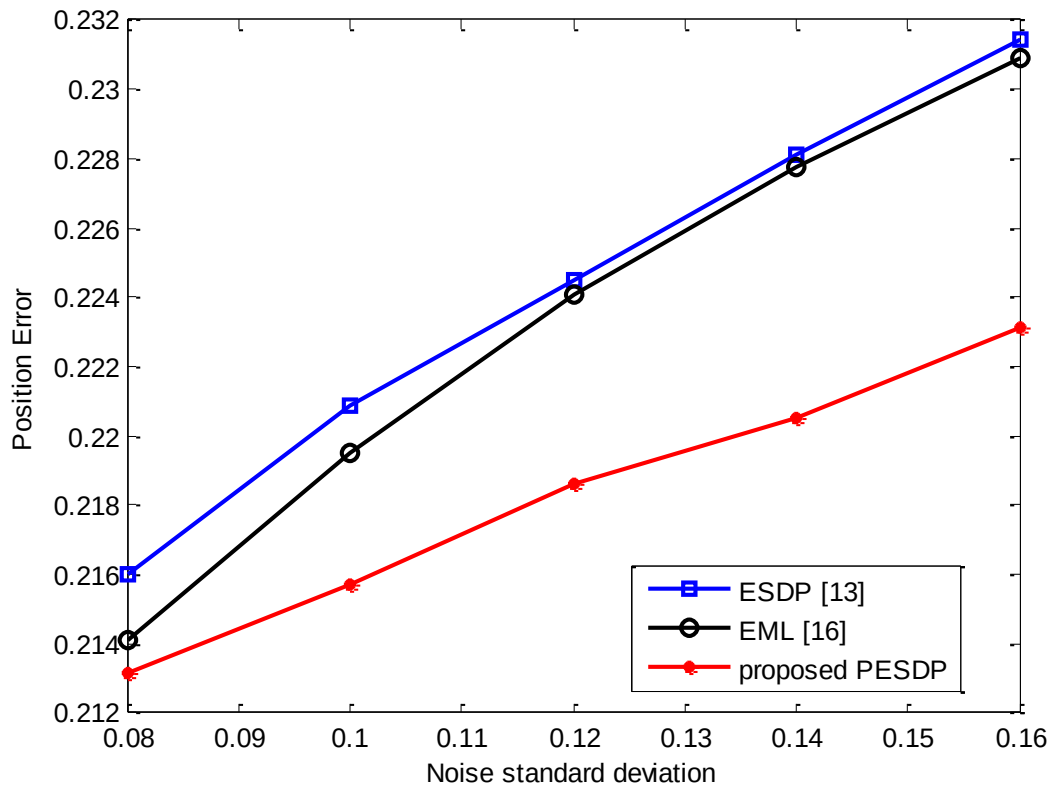

Fig. 1. Position Error of our proposed PESDP method (9), ESDP relaxation [13] and EML relaxation [16] in the presence of Gaussian noise.

## IV. CONCLUSIONS

In this paper, we proposed a less noise-sensitive convex relaxation (called PESDP) for wireless sensor network localization problem based on ESDP relaxation. PESDP provides a low rank solution to the problem as well as its dual. In PESDP, we modify ESDP model by perturbation matrix to make it less noise sensitive. By determination of an appropriate perturbation matrix, PESDP is compatible with all levels of noise. PESDP provides more accuracy in comparison with ESDP and EML especially when the noise level is high. Simulation results confirm that the complexity of the proposed PESDP is less than ESDP and EML methods, while it provides better accuracy.

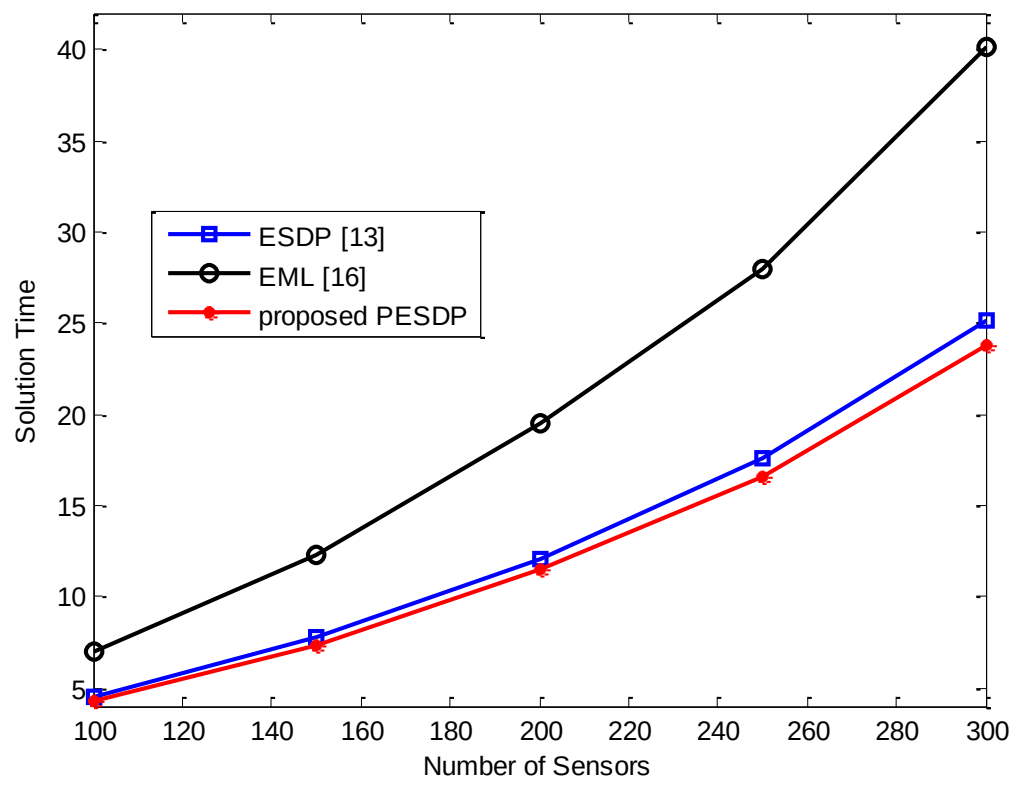

Fig. 2. Solution time of our proposed PESDP approach (9), ESDP relaxation [13] and EML relaxation [16].

## APPENDIX A

### PROOF OF (11)

Assume that $A \prec 0$. Therefore, by (7c) we may write the following statement:

$$S^{(i,j)}_{(1,2,i,j),(1,2,i,j)} - Z^T_{(1,2,i,j),(1,2,i,j)} A Z_{(1,2,i,j),(1,2,i,j)} \succcurlyeq 0 \tag{14}$$

Then, we may conclude (15) using Schur complement:

$$\begin{pmatrix} S^{(i,j)}_{(1,2,i,j),(1,2,i,j)} & Z^T_{(1,2,i,j),(1,2,i,j)} \\ Z_{(1,2,i,j),(1,2,i,j)} & A^{-1} \end{pmatrix} \succcurlyeq 0 \tag{15}$$

Now, consider the following lemma [19]:

Let $U \in \mathcal{R}^{m\times n}$ be a given matrix then $rank(U) \leq q$ if and only if there exist matrices $V = V^T \in \mathcal{R}^{m\times m}$ and $W = W^T \in \mathcal{R}^{n\times n}$ such that :

$$rank(V) + rank(W) \leq 2q$$

$$\begin{pmatrix} V & U \\ U^T & W \end{pmatrix} \succcurlyeq 0 \tag{16}$$

Therefore, we may conclude that:

$$rank\left(Z_{(1,2,i,j),(1,2,i,j)}\right) \leq q \Leftrightarrow$$

$$rank\left(S^{(i,j)}_{(1,2,i,j),(1,2,i,j)}\right) + 4 \leq 2q$$